\documentclass{article}

\usepackage{amsmath,amssymb,amsthm}
\usepackage{amsrefs}

\newtheorem{theorem}{Theorem}

\newtheorem{proposition}{Proposition}
\newtheorem{lemma}{Lemma}
\newtheorem{remark}{Remark}
\theoremstyle{definition}
\newtheorem{definition}{Definition}
\newtheorem{example}{Example}

\newcommand{\bkR}{\mathbb{R}}
\newcommand{\bkN}{\mathbb{N}}

\newcommand{\tH}{\mathcal{H}}
\newcommand{\tX}{\mathcal{X}}
\newcommand{\tT}{\mathcal{T}}
\newcommand{\bibx}[3]{}

\begin{document}

\title{Quadratures of Pontryagin Extremals\\
       for Optimal Control Problems\footnote{Presented at the
       4th Junior European Meeting on ``Control and Optimization'',
       Bia\l ystok Technical University, Bia\l ystok, Poland, 11-14 September 2005. Research Report CM05/I-48.}}

\author{Eug\'enio A. M. Rocha \and Delfim F. M. Torres}

\date{Department of Mathematics\\ University of Aveiro\\
        3810-193 Aveiro, Portugal\\
        \medskip
        \texttt{\{eugenio,delfim\}@mat.ua.pt}}

\maketitle

\begin{abstract}
We obtain a method to compute effective first integrals by
combining Noether's principle with the Kozlov-Kolesnikov
integrability theorem. A sufficient condition for the integrability by
quadratures of optimal control problems with controls taking
values on open sets is obtained. We illustrate our approach on
some problems taken from the literature. An alternative proof of
the integrability of the sub-Riemannian nilpotent Lie group of
type $(2,3,5)$ is also given.
\end{abstract}

\smallskip


\smallskip


\smallskip


\medskip


\section{Introduction}

Optimal control problems, with controls taking values on an open
set, are now subject to strong investigation because of their
recent applications to modern technology such as ``smart
materials'' \cite{MR2108960}. Geodesics of (sub)-Riemannian
manifolds can also be seen as solutions of a class of these
problems \cite{bonnard}. Meanwhile, solutions of optimal control
problems are closely related with solutions of Hamiltonian
equations through the Pontryagin Maximum Principle
\cite{MR0166037}. In the literature, Hamiltonian  equations are
usually classified as integrable or non-integrable. However, it was
not always clear in what sense a system is (non)integrable, as the
following quotations confirm \cite{dynamical,goriely}: Birkhoff comment
\emph{``when, however, one attempts to formulate a precise
definition of integrability, many possibilities appear, each with
a certain intrinsic theoretic interest''}; and the dictum of
Poincar\'e \emph{``a system of differential equations is only more
or less integrable''}. The reasons are the existence of
three main approaches to study integrability (the dynamical system approach
through bifurcation theory, the analytic approach using singularity analysis
and Painlev\'e property, and the algebraic approach through differential geometry),
which are not always compatible, and the stratification of the phase space
in regions that may be integrable with different notions.

On the Control Theory setting, integrability concerns the
existence of a foliation as the collection of all maximal integral
manifolds of the underlying distribution (see Frobenius and
Nagano-Sussmann theorems \cite{nagano}). On the other hand,
integration by quadratures of a differential equation is the
search for representing the solutions by a finite number of
algebraic operations, inversion of functions, and calculations of
integrals of known functions (``quadratures''); where a precise
meaning of integrability defines the allowed operations and the
set of known functions. Hence, an optimal control problem is integrable
by quadratures if the corresponding Hamiltonian equations are
integrable by quadratures. Since, an optimal control problem can
be integrable and not integrable by quadratures, and it is
frequent to shortcut integration by quadratures to integration, we
will adopt the word ``solvability'' to mean ``integrability by
quadratures'' of optimal control problems.

In the algebraic approach, verification of solvability of an optimal control
problem frequently requires the existence of a set of first integrals (or
conservation laws) of the true Hamiltonian equations, and an
appropriate method of reduction. The process can be divided in
three steps: (i) find a sufficient number of first integrals; (ii)
verify that such first integrals imply the existence of
quadratures; and (iii) apply a method to find the quadratures. The
first step can be attained by methods such as the Noether's
theorem, finding Casimirs, or even solving the PDE that appear in
the definition of variational symmetry. The second step can be
accomplished by using theorems such as Bour-Liouville,
Liouville-Arnold (abelian case), Mishchenko-Fomenko (nonabelian
case), or Kozlov-Kolesnikov. Last step is usually a consequence of
the choices on the first and second steps, \textrm{e.g.} Liouville-Poincar\'e
method, Cartan method, or Prykarpastsky method.
Notice that step (i), alone, is not enough to solve the problem,
since any $C^1$ function of a set of first integrals is a first
integral \cite[pp.~27]{goriely}, meaning that there exists plenty first
integrals that are useless.

In this work, we use the algebraic approach to derive in \S\ref{sec:mr:efi} a method for
computing effective first integrals (steps (i) and (ii)) for optimal control problems by
combining Noether's symmetry theorem with the Kozlov-Kolesnikov
integrability theorem (recalled in Section
\S\ref{sec:prelim}). Main result gives a sufficient
condition for the solvability of a given optimal control problem
(Theorem~\ref{main_teo}). A key issue is the construction of a
system of algebraic equations whose solutions determine the set of
effective first integrals, based on the simple observation that
for optimal control problems Noether's theorem usually gives a
parametric family of first integrals. The proposed method is
applied in \S\ref{sec:ie} to concrete optimal control problems
from the literature. An alternative proof to \cite{sachkov} for
the solvability of the sub-Riemannian nilpotent Lie group of type
$(2,3,5)$ is presented in \S\ref{sec:SR235}.


\section{Preliminaries}
\label{sec:prelim}

\subsection{The problem}

The optimal control problem consists to minimize a cost functional
\begin{equation}\label{eq:tr:prL:cs}
I\left[x(\cdot),u(\cdot)\right] = \int_a^b
\mathcal{L}\left(t,x(t),u(t)\right) dt
\end{equation}
subject to a control system described by ordinary differential equations
\begin{equation}
\label{eq:tr:prP:cs}
\dot{x}(t) = \varphi\left(t,x(t),u(t)\right)
\end{equation}
together with certain appropriate endpoint conditions. The
Lagrangian $\mathcal{L} : [a,b] \times \mathbb{R}^n \times
\mathbb{R}^m \rightarrow \mathbb{R}$ and the velocity vector
$\varphi : [a,b] \times \mathbb{R}^n \times \mathbb{R}^m
\rightarrow \mathbb{R}^n$ are given, and assumed to be smooth:
$\mathcal{L}(\cdot,\cdot,\cdot)$, $\varphi(\cdot,\cdot,\cdot)$
$\in C^1$. We are interested in the case where there the control
set is open: $u(t) \in U\subseteq\mathbb{R}^m$, with $U$ a open
set. We denote the problem by $(P)$. In the particular case
$\varphi(t,x,u) = u$, one obtains the fundamental problem of the
calculus of variations, which covers all classical mechanics. The
choice of the classes $\mathcal{X}$ and $\mathcal{U}$,
respectively of the state $x : [a,b] \rightarrow \mathbb{R}^n$ and
control variables $u : [a,b] \rightarrow \mathbb{R}^m$, are
important for the problem to be well-defined. We will assume, for
simplicity, that $\mathcal{X} =
PC^1\left([a,b];\mathbb{R}^n\right)$ and $\mathcal{U} =
PC\left([a,b];\mathbb{R}^m\right)$.

The Pontryagin Maximum Principle \cite{MR0166037} is a necessary
optimality condition which can be obtained from a general Lagrange
multiplier theorem in spaces of infinite dimension. Introducing
the Hamiltonian function
\begin{equation}
\label{eq:tr:H} H(x,u,\psi,t) = - \mathcal{L}(t,x,u) + \psi \cdot
\varphi(t,x,u) \, ,
\end{equation}
where $\psi_i$, $i = 1,\ldots,n$, are the ``Lagrange multipliers''
or the ``generalized momenta'', the multiplier theorem asserts that
the optimal control problem is equivalent to the maximization
of the augmented functional
\begin{equation*}
J\left[x(\cdot),u(\cdot),\psi(\cdot)\right] = \int_a^b \left(
H\left(x(t),u(t),\psi(t),t\right) - \psi(t) \cdot \dot{x}(t)
\right) dt \, .
\end{equation*}
Let $\left(\tilde{x}(\cdot),\tilde{u}(\cdot),\tilde{\psi}(\cdot)\right)$
solve the problem, and consider arbitrary $C^1$-functions
$h_1 \, , h_3 : [a,b] \rightarrow \mathbb{R}^n$,
$h_1(\cdot)$ vanishing at $a$ and $b$ ($h_1(\cdot) \in
C^1_0\left([a,b]\right)$),
and arbitrary continuous $h_2 : [a,b] \rightarrow \mathbb{R}^m$.
Let $\varepsilon$ be a scalar. By definition of maximizer, we have
\begin{equation*}
J\left[(\tilde{x}+\varepsilon h_1)(\cdot),(\tilde{u}+\varepsilon h_2)(\cdot),
(\tilde{\psi}+\varepsilon h_3)(\cdot)\right] \le
J\left[\tilde{x}(\cdot),\tilde{u}(\cdot),\tilde{\psi}(\cdot)\right] \, ,
\end{equation*}
and one has the following necessary condition:
\begin{equation}
\label{eq:tr:ncci}
\frac{d}{d\varepsilon} \left.
J\left[(\tilde{x}+\varepsilon h_1)(\cdot),(\tilde{u}+\varepsilon
h_2)(\cdot),(\tilde{\psi}
+\varepsilon h_3)(\cdot)\right]
\right|_{\varepsilon = 0} = 0\, .
\end{equation}
Differentiating \eqref{eq:tr:ncci} gives
\begin{multline*}
0 = \int_a^b \left[ \frac{\partial H}{\partial
x}\left(\tilde{x}(t), \tilde{u}(t),\tilde{\psi}(t),t\right) \cdot
h_1(t) + \frac{\partial H}{\partial u}\left(\tilde{x}(t),
\tilde{u}(t),\tilde{\psi}(t),t\right) \cdot h_2(t) \right. \\
\left. + \frac{\partial H}{\partial \psi}\left(\tilde{x}(t),
\tilde{u}(t),\tilde{\psi}(t),t\right) \cdot h_3(t)
- h_3(t) \cdot \dot{\tilde{x}}(t) - \tilde{\psi}(t) \cdot \dot{h}_1(t)\right] dt
\, .
\end{multline*}
Integrating the $\tilde{\psi}(t) \cdot \dot{h}_1(t)$ term by parts,
and having in mind that $h_1(a) = h_1(b) = 0$, one derives
\begin{multline}
\label{eq:tr:afterPP}
\int_a^b \left[ \left(\frac{\partial H}{\partial x}\left(\tilde{x}(t),
\tilde{u}(t),\tilde{\psi}(t),t\right) + \dot{\tilde{\psi}}(t)\right) \cdot
h_1(t) \right.
+ \frac{\partial H}{\partial u}\left(\tilde{x}(t),
\tilde{u}(t),\tilde{\psi}(t),t\right) \cdot h_2(t) \\
\left. + \left(\frac{\partial H}{\partial \psi}\left(\tilde{x}(t),
\tilde{u}(t),\tilde{\psi}(t),t\right) - \dot{\tilde{x}}(t)\right) \cdot
h_3(t)\right] dt = 0 \, .
\end{multline}
Note that \eqref{eq:tr:afterPP} was obtained for any variation $h_1(\cdot)$,
$h_2(\cdot)$, and $h_3(\cdot)$. Choosing $h_1(t) = h_2(t) \equiv 0$, and
$h_3(\cdot)$ arbitrary, one obtains the control system
\eqref{eq:tr:prP:cs}:
\begin{equation}
\label{eq:tr:PMP:cs}
\dot{\tilde{x}}(t) = \frac{\partial H}{\partial \psi}\left(\tilde{x}(t),
\tilde{u}(t),\tilde{\psi}(t),t\right) \, , \quad t \in [a,b] \, .
\end{equation}
With $h_1(\cdot)$ arbitrary, and $h_2(t) = h_3(t) \equiv 0$,
we obtain the \emph{adjoint system}:
\begin{equation}
\label{eq:tr:PMP:as}
\dot{\tilde{\psi}}(t) = - \frac{\partial H}{\partial x}\left(\tilde{x}(t),
\tilde{u}(t),\tilde{\psi}(t),t\right) \, , \quad t \in [a,b] \, .
\end{equation}
Finally, with $h_2(\cdot)$ arbitrary, and $h_1(t) = h_3(t) \equiv 0$,
the \emph{stationary condition} is obtained:
\begin{equation}
\label{eq:tr:PMP:mc}
\frac{\partial H}{\partial u}\left(\tilde{x}(t),
\tilde{u}(t),\tilde{\psi}(t),t\right) = 0 \, , \quad t \in [a,b] \, .
\end{equation}
Hence, a necessary optimality condition for
$\left(\tilde{x}(\cdot),\tilde{u}(\cdot)\right)$
to be a minimizer of problem $(P)$ is given by the \emph{Pontryagin Maximum
Principle}:
there exists $\tilde{\psi}(\cdot)$ such that the 3-tuple
$\left(\tilde{x}(\cdot),\tilde{u}(\cdot),\tilde{\psi}(\cdot)\right)$ satisfy
all the conditions \eqref{eq:tr:PMP:cs}, \eqref{eq:tr:PMP:as}, and
\eqref{eq:tr:PMP:mc}.
We recall that conditions \eqref{eq:tr:PMP:cs}, \eqref{eq:tr:PMP:as},
and \eqref{eq:tr:PMP:mc} imply the equality
\begin{equation}
\label{eq:tr:PMP:dHdt}
\frac{d}{dt} H\left(\tilde{x}(t),\tilde{u}(t),\tilde{\psi}(t),t\right)
= \frac{\partial H}{\partial
t}\left(\tilde{x}(t),\tilde{u}(t),\tilde{\psi}(t),t\right) \, .
\end{equation}

We assume, without loss of generality, that there exist at least
one $k\in\bkN_0$ such that
\begin{equation}
\exists\:(x,\psi,t)\in M : \frac{\partial}{\partial u}
\frac{d^{k}}{d t^{k}}\frac{\partial H}{\partial u}(x,u,\psi,t)
\not\equiv 0.
\end{equation}
Further, we denote by $k_u$ the smallest of such $k$ and by $\bar{u}$ the solution of
the equation
\begin{equation}
 \frac{d^{k_u}}{d t^{k_u}}\frac{\partial H}{\partial u}(x,u,\psi,t)=0
\end{equation}
with respect to $u$; giving as \emph{true Hamiltonian} the expression
\begin{equation}
 \tH(x,\psi,t)=H\left(x,\bar{u},\psi,t\right).
\end{equation}
A mapping $F(x,\psi,t)$ is a \emph{first integral} of the
Hamiltonian equations with Hamiltonian $\tH(x,\psi,t)$ if
\begin{equation}
 \frac{\partial F}{\partial t}+\{\tH,F\}=0,
\end{equation}
where $\{\cdot,\cdot\}$ denotes the canonical Poisson bracket.
\begin{remark}\label{rem01}
The restriction for the control set to be open is crucial. For
closed control sets $U$ the stationary condition
(\ref{eq:tr:PMP:mc}) become, the more general, maximality
condition $H(\tilde{x}(t),\tilde{u}(t),\tilde{\psi}(t))=\max_{v\in
U} H(\tilde{x}(t),v,\tilde{\psi}(t))$, and for such
cases the true Hamiltonian may be discontinuous. Our approach,
through symplectic geometry, can only deal with at least $C^1$
Hamiltonians.
\end{remark}

\begin{remark}\label{rem02}
If there exists a mapping $G(x,\psi)$ such that $\{\tH,G\}=c\tH$
for some $c\in\bkR$, then
$F(x,\psi,t)=G(x,\psi)-c\,t\,\tH(x,\psi)$ is a (nonautonomous)
first integral.
\end{remark}

\begin{remark}\label{rem03}
On this work we will consider nonautonomous problems, for such
reason we define the extended cotangent space of $\bkR^n$ by
$M=\bkR^n\{x\}\times(\bkR^n)^\ast\{\psi\}\times\bkR\{t\}$.
However, notice that a nonautonomous Hamiltonian on $\bkR^n$ is
not significantly different from an autonomous Hamiltonian on
$\bkR^{n+1}$, since a nonautonomous Hamiltonian $H(x,\psi,t)$ with
Hamiltonian system
$$   \dot{x}=\frac{\partial H}{\partial \psi}
\:\:\:\mbox{ and }\:\:\:
\dot{\psi}=-\frac{\partial H}{\partial x}
$$
can be transformed into an autonomous Hamiltonian
$K(x,\psi,\theta,t)=H(x,\psi,t)-\theta$
with Hamiltonian equations
$$   \dot{x}=\frac{\partial K}{\partial \psi},\:\:\:
\dot{\psi}=-\frac{\partial K}{\partial x},\:\:\:
\dot{\theta}=\frac{\partial K}{\partial t}
\:\:\:\mbox{ and }\:\:\:
\dot{t}=-\frac{\partial K}{\partial \theta}.
$$
The reverse procedure partially justifies the statement that an
autonomous Hamiltonian is always a first integral for the problem,
so the equations of an autonomous Hamiltonian can be
dimension-reduced. Furthermore, a system with only one degree of
freedom is always integrable.
\end{remark}


\subsection{Noether's theorem}

In 1918 Emmy Noether established the key result to find conservation
laws in the calculus of variations \cite{MR0406752}. We sketch here
the standard argument used to derive Noether's theorem and conservation laws
in the optimal control setting (\textrm{cf.} \textrm{e.g.}
\cite{MR1901565,MR0341229}).

Let us consider a one-parameter group of $C^1$-transformations of the form
\begin{equation}
\label{eq:tr:pt} h_s(x,u,\psi,t) =
\left(h^x_s(x,u,\psi,t),h^u_s(x,u,\psi,t),h^\psi_s(x,u,\psi,t),h^t_s(x,u,\psi,t)\right)
\, ,
\end{equation}
where $s$ denote the independent parameter of the transformations.
We require that to the parameter value $s = 0$ corresponds
the identity transformation:
\begin{equation}
\label{eq:tr:pts0}
\begin{split}
h_0(x,u,\psi,t) &= \left(h^x_0(x,u,\psi,t),h^u_0(x,u,\psi,t),
h^\psi_0(x,u,\psi,t),h^t_0(x,u,\psi,t)\right) \\
&= (x,u,\psi,t) \, .
\end{split}
\end{equation}
Associated to the group of transformations \eqref{eq:tr:pt} we consider
the \emph{infinitesimal generators}
\begin{gather}
T(x,u,\psi,t) = \left.\frac{d}{ds} h^t_s(x,u,\psi,t)\right|_{s = 0} \, , \quad
X(x,u,\psi,t) = \left.\frac{d}{ds} h^x_s(x,u,\psi,t)\right|_{s = 0} \, ,
\label{eq:tr:TX} \\
U(x,u,\psi,t) = \left.\frac{d}{ds} h^u_s(x,u,\psi,t)\right|_{s = 0} \, , \quad
\Psi(x,u,t,\psi) = \left.\frac{d}{ds} h^\psi_s(x,u,\psi,t)\right|_{s = 0} \, .
\notag
\end{gather}

\begin{definition}
\label{def:tr:inv}
The optimal control problem $(P)$ is said to be invariant under
a one-parameter group of $C^1$-transformations \eqref{eq:tr:pt}
if, and only if,
\begin{multline}
\label{eq:tr:inv}
\left.\frac{d}{ds}
  \Biggl\{
    \biggl[
      H\left(h_s\left(x(t),u(t),\psi(t),t\right)\right) \right. \\
      \left. - h^\psi_s\left(x(t),u(t),\psi(t),t\right) \cdot
\frac{\frac{dh^x_s\left(x(t),u(t),\psi(t),t\right)}{dt}}{\frac{dh^t_s\left(x(t),
u(t),\psi(t),t\right)}{dt}}
    \biggr]
    \frac{dh^t_s\left(x(t),u(t),\psi(t),t\right)}{dt}
  \Biggr\}
\right|_{s = 0} = 0 \, ,
\end{multline}
with $H$ the Hamiltonian \eqref{eq:tr:H}.
\end{definition}

Having in mind \eqref{eq:tr:pts0}, condition \eqref{eq:tr:inv}
is equivalent to
\begin{multline}
\label{eq:tr:cond:nec:suf:inv}
\frac{\partial H}{\partial t} T +
\frac{\partial H}{\partial x} \cdot X + \frac{\partial H}{\partial u} \cdot U +
\frac{\partial H}{\partial \psi} \cdot \Psi - \Psi \cdot \dot{x}(t)
- \psi(t) \cdot \frac{d}{dt}X + H  \frac{d}{dt} T = 0 \, ,
\end{multline}
where here all functions are evaluated
at $\left(x(t),u(t),\psi(t),t\right)$ whenever not indicated.
Along a Pontryagin extremal $\left(x(\cdot),u(\cdot),\psi(\cdot)\right)$
equalities \eqref{eq:tr:PMP:cs}, \eqref{eq:tr:PMP:as}, \eqref{eq:tr:PMP:mc},
and \eqref{eq:tr:PMP:dHdt} are in force, and \eqref{eq:tr:cond:nec:suf:inv}
reduces to
\begin{equation}
\label{eq:ExtPont:cond:nec:suf:inv}
\frac{dH}{dt} T -
\dot{\psi}(t) \cdot X - \psi(t) \cdot \frac{dX}{dt} + H
\frac{dT}{dt} = 0 \Leftrightarrow \frac{d}{dt} \left( \psi(t)
\cdot X - H T\right) = 0 \, .
\end{equation}
Therefore, we have just proved \emph{Noether's theorem for optimal control problems}.
\begin{theorem}[Noether's Theorem]\label{noether_theo}
If the optimal control problem is invariant under \eqref{eq:tr:pt},
in the sense of Definition~\ref{def:tr:inv}, then
\begin{equation}
\label{eq:tr:cl}
\psi(t) \cdot X\left(x(t),u(t),\psi(t),t\right)
- H\left(x(t),u(t),\psi(t),t\right) T\left(x(t),u(t),\psi(t),t\right) = const
\end{equation}
($t \in [a,b]$; $T$ and $X$ are given according to
\eqref{eq:tr:TX};
$H$ is the Hamiltonian \eqref{eq:tr:H}) is a \emph{conservation law},
that is, \eqref{eq:tr:cl} is valid along all the minimizers
$(x(\cdot),u(\cdot))$ of $(P)$ which are Pontryagin extremals.
\end{theorem}


\subsection{Solvability and Reduction}

E.~Bour and J.~Liouville, on the middle of the XIX century,
obtained fundamental concepts and results concerning integrability
by quadratures of differential equations. Namely, the notion of
elementary function: a $C^n$ function that belong to the set
$\Lambda$ of elementary functions. The set $\Lambda$ is obtained
from rational functions on $C^k$ ($k\in\bkN_0$), using a finite
number of the following operations: (i) algebraic operations (if
$f_1,f_2\in\Lambda$ then $f_1\star f_2\in\Lambda$, where $\star$
is either the addition, subtraction, multiplication, or
division); (ii) solutions of algebraic equations with coefficients
in $\Lambda$; (iii) differentiation; and (iv) exponential and
logarithm operations. The set of elementary functions together with
the operation of integration (if $f\in\Lambda$ then $\int
f(x)\,dx\in\Lambda$) is called the set of Liouvillian functions.
Liouville (1939) showed that the solution of the equation
$\dot{x}(t)=t^\alpha-x^2$ is only Liouvillian for $\alpha=-2$ and
$\alpha=4k/(1-2k)$ ($k\in\bkN$).

There is a concrete method that permits not only to verify that
solutions of the Hamiltonian system are Liouvillian functions,
but also to reduce the system in order to obtain the extremals. We
describe it briefly. Let
$M=\bkR^n\{x\}\times(\bkR^n)^\ast\{\psi\}\times\bkR\{t\}$ and
$f:M\rightarrow\bkR$ be a first integral of the Hamiltonian system
with Hamiltonian $\tH$. If $df(q)\neq 0$, then in some
neighborhood of the point $q\in M$ there exist symplectic
coordinates $(\tilde{x},\tilde{\psi},\tilde{\theta},\tilde{t})$
such that $f(\tilde{x},\tilde{\psi}, \tilde{t})=\tilde{\psi}_1$.
In these coordinates $\tH$ does not depend on $\tilde{x}_1$,
therefore if we fix a value $f=\tilde{\psi}_1=c$, then the
Hamiltonian system will only have $n-1$ degrees of freedom.
In order to have an
effective reduction of dimension by a set of first integrals, this
method requires the first integrals to be independent and in
involution
\begin{equation*}
 \{f_i,f_j\}=0\:\:\: \forall\,i,j\in\{1,\dots,N\}.
\end{equation*}

E.~Cartan \cite{Cartan} extended Liouville's method for the case
where the algebra $L$ of first integrals is not commutative, and
possible infinite-dimensional. Cartan assumed that the first
integrals satisfy the relation
\begin{equation}\label{CartanCond}
  \{f_i,f_j\}=\zeta_{ij}(f_1,\dots,f_N)\:\:\:\forall\, i,j\in\{1,\dots,N\}
\end{equation}
for some (nonlinear) functions $\zeta_{ij}:\bkR^N\rightarrow\bkR$.
His method is based on the following theorem.

\begin{theorem}[S.Lie - E.Cartan] \cite{Cartan}
\label{th:cartan} Let $F=(f_1,\dots,f_N)$. Suppose that the point
$c\in\bkR^N$ is not a critical value of the mapping $F$ and that
in its neighborhood the rank of the matrix $(\zeta_{ij})$ is
constant. Then, in a small neighborhood $U\subset\bkR^N$ of $c$,
one can find $N$ independent functions $\phi_j:U\rightarrow\bkR$
such that the functions $\Phi_j=\phi_j\circ F: V\rightarrow\bkR$,
where $V=F^{-1}(U)$, satisfy the relations
$$
 \{\Phi_1,\Phi_2\}=\dots=\{\Phi_{2\eta-1},\Phi_{2\eta}\}=1,
$$
whereas the remaining brackets vanish, and the rank of the matrix
$(\zeta_{ij})$ is $2\eta$.
\end{theorem}
Using Theorem~\ref{th:cartan} we can lower the order of the system
in the following way: the level set $M_c=\{(x,\psi)\in
M:\Phi_j(x,\psi)=c_j, 1\leq j\leq N\}$, where $c=(c_1,\dots,c_N)$
satisfies the theorem, is a smooth $(2n-N)$-dimensional
submanifold of $M$. The theorem also implies that there is an
action of the commutative group $\bkR^l$ ($l=N-2\eta$) on $M_c$,
generated by the phase flows of the Hamilton's equations with
Hamiltonians $\Phi_j$ for $j>2\eta$. Now, thanks to the functional
independence of the integrals $\Phi_j$, this action has no fixed
points. Hence, if its orbits are compact, then the quotient space
$M_{red}=M_c/\bkR^l$ is a smooth manifold with dimension
$2(n-N+\eta)$ endowed with a natural symplectic structure. Let
$H'$ denote the restriction of the Hamiltonian $H$ to the level
set $M_c$ of the first integrals. Since $H'$ is constant on the
orbits of the group $\bkR^l$, there is a smooth function
$H_{red}:M_{red}\rightarrow\bkR$ such that the diagram
$$
M_c\stackrel{pr}{\longrightarrow}M_{red}\stackrel{H_{red}}{\longrightarrow}\bkR
 \stackrel{H'}{\longleftarrow}M_c$$
commutes. To end, let us observe that one can obtain Liouville's
method from Cartan's method by choosing $\zeta_{ij}\equiv 0$ and
$\eta=0$. Locally they give the same result, however the
factorization by Cartan method can be accomplished globally only
under more restrictive assumptions.

Kozlov-Kolesnikov \cite{KK} proved an intermediate result,
considering that Poisson brackets of first integrals are a linear
combination of first integrals. In fact, this result is more suitable for
our purposes.

\begin{theorem}[Kozlov-Kolesnikov]
Suppose that the Hamiltonian
$\mathcal{H}:\bkR^{n}\times(\bkR^{n})^\ast\times\bkR\rightarrow\bkR$ has $n$
first integrals
$F_1,\dots,F_n:\bkR^{n}\times(\bkR^{n})^\ast\times\bkR\rightarrow\bkR$ that
satisfy relation
\begin{equation}\label{KozlovInv}
  \exists\, \xi^{ij}\in\bkR^n \:\: : \:\:\{F_i,F_j\}=\sum_{s=1}^n \xi^{ij}_s
F_s\:\:\:\:\forall i,j\in\{1,\dots,n\},
\end{equation}
where $\xi^{ij}=(\xi^{ij}_1,\dots,\xi^{ij}_n)^T$.
Additionally, assume that
\begin{enumerate}
\item on the set $M_f=\{(x,\psi,t)\in M: F_i(x,\psi,t)=r_i, 1\leq i\leq n\}$ the functions
$F_1,\dots,F_n$ are independent;
\item $(r_1,\dots,r_n)\xi^{ij}=0$ for all $i,j=\{1,\dots,n\}$;
\item the Lie algebra $L$ of linear combination $\sum_i c_iF_i$, $c_i\in\bkR$,
is solvable.
\end{enumerate}
Then the solutions of the Hamiltonian system lying on $M_f$ and can be found by
quadratures.
\end{theorem}


\section{Main Results -- effective first integrals}
\label{sec:mr:efi}

In order to use Noether's theorem (Theorem~\ref{noether_theo})
to obtain effective first
integrals for the Hamiltonian equations, we need to compute the
solutions of the first order partial differential equation
(\ref{eq:ExtPont:cond:nec:suf:inv}) considering the optimal
control $\bar{u}$, \textrm{i.e.}
\begin{multline}
\label{eq:tinv} \frac{d\tH(x,\psi,t)}{dt} \tT(x,\psi,t) -
\dot{\psi}(t) \cdot \tX(x,\psi,t)\\ - \psi(t) \cdot
\frac{d\tX(x,\psi,t)}{dt} + \tH(x,\psi,t)
\frac{d\tT(x,\psi,t)}{dt} = 0,
\end{multline}
where $\tX(x,\psi,t)\equiv X\left(x,\bar{u},\psi,t\right)$ and
$\tT(x,\psi,t)\equiv T\left(x,\bar{u},\psi,t\right)$.
A particular solution (if exists) can be found, \textrm{e.g.}
by the well known method of (additive) separation of variables for PDE, assuming
\begin{eqnarray*}
  \tT(x,\psi,t) &=& \tT^0(t)+\tT^{x_1}(x_1)+\dots
  +\tT^{x_n}(x_n)+\tT^{\psi_1}(\psi_1)+\dots+\tT^{\psi_n}(\psi_n),\\
  \tX(x,\psi,t) &=& \tX^0(t)+\tX^{x_1}(x_1)+\dots
  +\tX^{x_n}(x_n)+\tX^{\psi_1}(\psi_1)+\dots+\tX^{\psi_n}(\psi_n).
\end{eqnarray*}
In particular, we will consider that each independent component of $\tT$ and $\tX$
have a polynomial structure with degree $\leq p_d$, hence
\begin{eqnarray}\label{solform:I}
  \tT(x,\psi,t) &=& \sum_{\mu=0}^{p_d}\left(C^{\tT}_0(\mu)t^\mu
   +\sum_\nu C^{\tT}_\nu(\mu) (q_\nu)^\mu\right),\\
  \tX_i(x,\psi,t) &=& \sum_{\mu=0}^{p_d}\left(C^{\tX}_0(i,\mu)t^\mu
   +\sum_\nu C^{\tX}_\nu(i,\mu) (q_\nu)^\mu\right),\label{solform:II}
\end{eqnarray}
for some constants $C^{\tT}_0(\mu),C^{\tT}_\nu(\mu),C^{\tX}_0(i,\mu),
C^{\tX}_\nu(i,\mu)\in\bkR$ and $\tX=(\tX_1,\dots,\tX_n)$.

Now, since equation~(\ref{eq:tinv}) have to be valid for every
extremal $(x(t),\psi(t),t)\in M$, this particular choice for the
structure of solutions will transform the PDE problem for
$\tT,\tX$ into an algebraic system of equations for the constants
$C^{\tT}_0(\mu)$, $C^{\tT}_\nu(\mu)$, $C^{\tX}_0(i,\mu)$,
$C^{\tX}_\nu(i,\mu)\in\bkR$. The algebraic system is
under-determined, because we have one equation for two unknowns
$\tT$ and $\tX$. Therefore, if such system has a nontrivial
solution, we have a family of first integrals. Let us define the
mapping $\digamma:\bkR^m\rightarrow C(M,\bkR)$ by
\begin{equation}
 \digamma(\lambda)(x,\psi,t)=\psi\cdot \tX - \tH \tT
\end{equation}
for $\lambda\in\bkR^m$. $\digamma$ is a linear mapping, whereas
$\tT$ and $\tX$ are linear with respect to the constants, the PDE
(\ref{eq:tinv}) is a linear first order equation (superposition of
solutions is a solution), and equation~(\ref{eq:tr:cl}) is a
linear functional combination of $\tT$ and $\tX$. We resume our
statements in the following lemma.

\begin{lemma}\label{lem01}
Assume that equation~(\ref{eq:tinv}) has nontrivial solutions
$\tT(x,\psi,t)$ and $\tX(x,\psi,t)$, respectively of the form
(\ref{solform:I}) and (\ref{solform:II}). Then the mapping
$\digamma$ is a linear $m$-parametric family of first integrals
with respect to $\lambda\in\bkR^m$, \textrm{i.e.} it depends on
$m\in\bkN_0$ arbitrary constants $\lambda_1,\dots,\lambda_m$ and
\begin{equation}
 \forall \lambda\in\bkR^m :  \frac{\partial \digamma(\lambda)}{\partial t}
  +\{\tH, \digamma(\lambda)\}=0.
\end{equation}
\end{lemma}

Considering the previous lemma, makes sense to assume
that $\digamma$ has the following structure
\begin{equation}\label{eq:digamma}
  \digamma(\lambda)(x,\psi,t)=\sum_{k=1}^m \digamma_k(x,\psi,t)\lambda_k,
\end{equation}
where $\lambda=(\lambda_1,\dots,\lambda_m)\in\bkR^m$. We call
$\digamma_k$ the components of the family of first integrals
$\digamma(\lambda)$.

Let $n$ be the dimension of the phase space ($x\in\bkR^n$),
$m$ the number of parameters on the family $\digamma(\lambda)$,
and $m\geq n$. The existence of $n$ effective first integrals
(in the sense of Kozlov-Kolesnikov) will be related with the existence of a
nontrivial solution of a system of algebraic equations involving
the components of $\digamma$ and their canonical Poisson brackets.
Consider the following system of algebraic equations
\begin{equation}\label{prob}
   (\lambda^i)^T\: A(x,\psi,t)\: \lambda^j = (\xi^{ij})^T
       \left[\begin{array}{c} (\lambda^1)^T \\ \vdots \\ (\lambda^n)^T
       \end{array}\right] b(x,\psi,t),\:\:\:\:
   i<j\in\{1,\dots,n\},
\end{equation}
where $\lambda^1,\dots,\lambda^n\in\bkR^m$ and $\xi^{ij}\in\bkR^n$,
$$
 A(x,\psi,t) = \left(\left\{\digamma_p,\digamma_q\right\}\right)_{p,q=1}^m
 \:\:\:\mbox{ and }\:\:\:
 b(x,\psi,t) = \left[ \digamma_1,\dots, \digamma_m\right]^T.
$$
By a solution of the system of equations~(\ref{prob}) we mean a
set of constant vectors $(\lambda^1,\dots,\lambda^n)$,
$(\xi^{12},\dots,\xi^{(n-1)n})$ that satisfies the system. A
nontrivial solution is a solution for which all $\lambda^k$
($k=1,\dots,n$) are different.
From a nontrivial solution,  we have the set of first integrals
$\{\digamma(\lambda^1),\dots,\digamma(\lambda^n)\}$. Comparing with
relation~(\ref{CartanCond}), $\zeta_{ij}$ are linear functions.
Therefore, the space of linear combinations
$L=span\{\digamma(\lambda^1),\dots,\digamma(\lambda^n)\}$
forms a noncommutative but finite-dimensional Lie algebra, where
the first integrals define a basis and the coordinates of $\xi^{ij}$
are the structure constants.

\begin{proposition}
\label{prop:1} If there exists a nontrivial solution
$(\lambda^1,\dots,\lambda^n)$, $(\xi^{12},\dots,\xi^{(n-1)n})$ to
the system~(\ref{prob}), then the set of first integrals
$\{\digamma(\lambda^1),\dots,\digamma(\lambda^n)\}$ satisfy
relation~(\ref{KozlovInv}).
\end{proposition}
\begin{proof}
The set $\{\digamma(\lambda^1),\dots,\digamma(\lambda^n)\}$ is a
set of first integrals of the Hamiltonian $\tH$, since
$\digamma(c_1,\dots,c_m)(x,\psi,t)$ is a first integral for any
$(c_1,\dots,c_m)\in\bkR^m$; by Theorem~\ref{noether_theo} and
Lemma~\ref{lem01}. Relation~(\ref{KozlovInv}) is satisfied, using
the definition of $\digamma$ (\ref{eq:digamma}), properties of the
bracket, and relation~(\ref{prob}). For
$i<j\in\{1,\dots,n\}$, we have: {\small
\begin{multline*}
\{\digamma(\lambda^i),\digamma(\lambda^j)\} =
  \left\{\sum_{p=1}^m \digamma_p\lambda^i_p,\sum_{q=1}^m
\digamma_q\lambda^j_q\right\}
   = \sum_{p=1}^m\sum_{q=1}^m
   \lambda^i_p\lambda^j_q\{\digamma_p,\digamma_q\}
= (\lambda^i)^T A(x,\psi,t) \lambda^j \\
 = (\xi^{ij})^T
       \left[\begin{array}{c} (\lambda^1)^T \\ \vdots \\ (\lambda^m)^T
       \end{array}\right] b(x,\psi,t)
 = (\xi^{ij})^T
       \left[\begin{array}{c} \digamma(\lambda^1) \\ \vdots \\
\digamma(\lambda^m)
       \end{array}\right] = \sum_{s=1}^n \xi^{ij}_s
       \digamma(\lambda^s).
\end{multline*}}
\end{proof}

\begin{proposition}
\label{prop:2} Assume that the set of first integrals
$\{\digamma(\lambda^1),\dots,\digamma(\lambda^n)\}$ satisfy
relation~(\ref{KozlovInv}). Let
$\mathcal{S}=\{(a,b,p,q,i,j)\in\{1,\dots,n\}^4\times\{1\dots,m\}^2:
 a<b, p<q, a<p, i<j \}$. If
\begin{equation}
 \xi^{ab}_i\xi^{pq}_j=\xi^{pq}_i\xi^{ab}_j\:\:\:\:\forall\:(a,b,p,q,i,j)\in
 \mathcal{S},
\end{equation}
then the Lie algebra $L$ of linear combination $\sum_s
c_s\digamma(\lambda^s)$, $c_s\in\bkR$, is solvable.
\end{proposition}
\begin{proof}
 Let $L^0\equiv L$.
 We recall that a Lie algebra $L$ is solvable if the descent
 series is nilpotent, i.e.
 \begin{equation}
   \exists \bar{k}\in\bkN : L^{\bar{k}}\equiv 0
   \:\:\:\mbox{ where }\:\:\:L^{k}=[L^{k-1},L^{k-1}].
 \end{equation}
 For our purpose it will be enough to consider $\bar{k}=2$. Let us observe that the
 Liouville method (first integrals in involution) is the case $\bar{k}=1$.
Hence, for $k=1$ and using relation~(\ref{KozlovInv}), {\small
 \begin{eqnarray*}
  \left\{\sum_a \alpha_a\digamma(\lambda^a),\sum_b \beta_b\digamma(\lambda^b)\right\}
  &=& \sum_{a<b}
  (\alpha_a \beta_b-\beta_a \alpha_b)\{\digamma(\lambda^a),\digamma(\lambda^b)\}\\
  &=& \sum_i\sum_{a<b}
  (\alpha_a \beta_b-\beta_a \alpha_b)\xi^{ab}_i\digamma(\lambda^i),
 \end{eqnarray*}}
and, for $k=2$, {\small
 \begin{eqnarray*}
  &\ & \left\{\sum_i\sum_{a<b}
  (\alpha_a \beta_b-\beta_a \alpha_b)\xi^{ab}_i\digamma(\lambda^i),
  \sum_j\sum_{p<q}
  (\alpha_p \beta_q-\beta_p \alpha_q)\xi^{pq}_j\digamma(\lambda^j)\right\}\\
  &\ & = \sum_s\sum_{i<j}
  \left(\left(\sum_{a<b}
  (\alpha_a \beta_b-\beta_a \alpha_b)\xi^{ab}_i\right)\left(\sum_{p<q}
  (\alpha_p \beta_q-\beta_p \alpha_q)\xi^{pq}_j\right)\right.\\
  &\ &\left.-\left(\sum_{p<q}
  (\alpha_p \beta_q-\beta_p \alpha_q)\xi^{pq}_i\right) \left(\sum_{a<b}
  (\alpha_a \beta_b-\beta_a
  \alpha_b)\xi^{ab}_j\right)\right)\xi^{ij}_s\digamma(\lambda^s)\\
  &\ & = \sum_s\sum_{i<j}
  \sum_{a<b}\sum_{p<q}(\alpha_a \beta_b-\beta_a \alpha_b)(\alpha_p \beta_q-\beta_p \alpha_q)
\left(\xi^{ab}_i\xi^{pq}_j-\xi^{pq}_i\xi^{ab}_j\right)\xi^{ij}_s\digamma(\lambda^s).
\end{eqnarray*}}
\end{proof}

We are now in conditions to present the main result of the paper:
a practical method to find effective first integrals for optimal
control problems. Theorem~\ref{main_teo} is a direct consequence
of Propositions~\ref{prop:1} and \ref{prop:2}, and
Kozlov-Kolesnikov theorem.
\begin{theorem}
\label{main_teo} Assume that the optimal control problem
(\ref{eq:tr:prL:cs})-(\ref{eq:tr:prP:cs}) has a solution,
and there exists a $m$-parametric family
of first integrals $\digamma$, given by Lemma~\ref{lem01}, with
the form
\begin{equation*}
  \digamma(\lambda)(x,\psi,t)=\sum_{k=1}^m \digamma_k(x,\psi,t)\lambda_k,
\end{equation*}
where $\lambda=(\lambda_1,\dots,\lambda_m)\in\bkR^m$.
Let $\mathcal{S}=\{(a,b,p,q,i,j)\in\{1,\dots,n\}^4\times\{1\dots,m\}^2:
 a<b, p<q, a<p, i<j \}$,  $A(x,\psi,t) =
 \left(\left\{\digamma_p,\digamma_q\right\}\right)_{p,q=1}^m$,
 $\Lambda=\left[(\lambda^1)^T,\dots,(\lambda^n)^T\right]^T\in M_{n\times m}$,
 and $b(x,\psi,t) = \left[ \digamma_1,\dots, \digamma_m\right]^T$.

\textbf{If} there exists a solution
$(\lambda^1,\dots,\lambda^n,\xi^{12},\dots,\xi^{(n-1)n},r_1,\dots,r_n)$ (with
$\lambda^i\in\bkR^m$, $\xi^{ij}\in\bkR^n$, $r_i\in\bkR$ and
$i<j\in\{1,\dots,n\}$) of the
algebraic system of equations
$$
\left\{\begin{array}{rcll}
(\lambda^i)^T\: A(x,\psi,t)\: \lambda^j - (\xi^{ij})^T
       \Lambda\, b(x,\psi,t) & = & 0,
       &\mbox{ for } \forall\,i<j\in\{1,\dots,n\},\\
 \xi^{ab}_i\xi^{pq}_j -\xi^{pq}_i\xi^{ab}_j & = &
 0, &\mbox{ for }\forall\:(a,b,p,q,i,j)\in \mathcal{S},\\
 \sum_{s=1}^n r_s\xi^{ij}_s & = & 0,
  &\mbox{ for } \forall\,i<j\in\{1,\dots,n\},
\end{array}\right.
$$
and
\begin{equation}
rank\,[\nabla_{(x,\psi)}\digamma(\lambda^1),\dots,\nabla_{(x,\psi)}\digamma(\lambda^n)]=n
\end{equation}
on the manifold
$M_\digamma=\{\alpha\in M: \digamma(\lambda^i)(\alpha)=r_i\},$
\textbf{then} the optimal control problem
(\ref{eq:tr:prL:cs})-(\ref{eq:tr:prP:cs}) is solvable on $M_\digamma$.
\end{theorem}

Although simple, the arguments behind Lemma~\ref{lem01} and
Theorem~\ref{main_teo} give a powerful method that can be applied
with success to several problems of optimal control.


\section{Illustrative Examples}
\label{sec:ie}

We now present three interesting applications, many others can
be chosen from the literature. Families of first integrals
(Lemma~\ref{lem01}) were obtained using the \textit{Maple}
package described in \cite{GouveiaTorresRocha}.


\begin{example}
Let us show the integrability by quadratures (solvability) of the
following optimal control problem
\begin{eqnarray*}
\frac{1}{2} \int_a^b u_1(t)^2+u_2(t)^2\,dt \rightarrow \min \, , \quad
\left\{\begin{array}{l}
\dot x_1(t)=u_1(t) \cos (x_3(t))\\
\dot x_2(t)=u_1(t) \sin (x_3(t))\\
\dot x_3(t)=u_2(t)
\end{array}\right.
\end{eqnarray*}
which is known as the Dubin's model for the kinematics of a car
\cite[Example~18, pp.~750--751]{Rouchon01}. The true Hamiltonian is
$$\tH=\frac{1}{2}\left(\left[\cos \left(
x_{3}(t)  \right) \psi_{1}(t)+\sin \left( x_{3}(t)\right)
\psi_{2}(t)\right]^2 +\left[\psi_{3} (t)\right]^2\right).$$ It is
clear that the problem has three trivial first integrals (f.i.)
$\{\tH,\psi_1,\psi_2\}$ in involution. Notice that an autonomous
Hamiltonian is always a f.i. by Remark~\ref{rem03}, and the other
f.i. follow from Remark~\ref{rem02}. However, by computing
$\digamma(\lambda)$ on Lemma~\ref{lem01} and applying
Theorem~\ref{main_teo}, we obtain the trivial f.i. and an extra
effective f.i. $F=-\psi_1 x_2+\psi_2 x_1 + \psi_3$. Therefore,
the set $\{\psi_1,\psi_2,F\}$ can also be
used to prove the solvability of the problem.
\end{example}


\begin{example}
An interesting variation of the previous problem is the model of a
car with one-trailer \cite{Fuka}, parameterized by
constants $(a,b,c)\in\bkR$,
\begin{eqnarray*}
\int_a^b u_1^2+u_2^2\,dt \rightarrow \min \, , \quad
\left\{\begin{array}{l}
\dot x_1 = u_1 \cos (x_3)\\[0.1cm]
\dot x_2 = u_1 \sin (x_3)\\[0.1cm]
\dot x_3 = \frac{1}{c} u_1\tan(u_2)\\[0.1cm]
\dot x_4 = \frac{1}{b}u_1\left(\frac{a}{c}\tan(u_2)\cos(x_3-x_4)
             -\sin(x_3-x_4)\right)
\end{array}\right.
\end{eqnarray*}
The necessary and sufficient condition of invariance is satisfied
with the following generators $\left\{ \tT=C_{2},
\tX_{1}=-C_{1}x_{2}+C_{4},
 \tX_{2}=C_{1}x_{1}+C_{3}, \tX_{3}=C_{1},\tX_{4}=C_{1}\right\}$.
 It follows that, for
$C=(C_1,C_2,C_3,C_4)\in\bkR^4$,
$$\digamma(C)(x,\psi,t) =
\left(C_4-C_1 x_2\right)\psi_1+ \left(C_3+C_1
x_1\right)\psi_2+C_1\psi_3(t)+C_1\psi_4 -C_{2}\tH.$$ Therefore, a
possible solution of the algebraic system of
Theorem~\ref{main_teo} is $\lambda^1=(1,0,0,0)$,
$\lambda^2=(0,0,1,0)$, $\lambda^3=(0,0,0,1)$,
$\lambda^4=(0,1,0,0)$, $\xi^{12}=(0,0,-1,0)$,
$\xi^{13}=(0,1,0,0)$, $\xi^{23}=(0,0,0,0)$, $\xi^{i4}=(0,0,0,0)$
for $i=1,2,3$, and $r=(c_1,0,0,c_4)$ for any $c_1,c_4\in\bkR$. The
set of effective first integrals is
$$\{\digamma(\lambda^1)=-\psi_1 x_2+\psi_2
x_1+\psi_3+\psi_4,\digamma(\lambda^2)=\psi_2,
\digamma(\lambda^3)=\psi_1,\digamma(\lambda^4)=\tH\}.$$ Comparing
with the last example, which is linear in the control, this
problem is not only nonlinear in the control as it has a very
complicated true Hamiltonian. However, the set of first integrals
is just an extension of the previous one, where the only change is
$\digamma(\lambda^1)=F+\psi_4$.
\end{example}


\begin{example}
We now consider the so-called flat Martinet problem
\cite{Trelat98}:
\begin{eqnarray*}
\int_a^b u_1^2+u_2^2\,dt \rightarrow \min \, , \quad
\left\{
\begin{array}{l}
\smallskip
\dot x_1=u_1\\
\displaystyle \smallskip \dot x_2=\frac{u_2}{1+\alpha x_1} \, ,\\
\dot x_3=x_2^2 u_1
\end{array}
\right.
\qquad \alpha \in \mathbb{R} \, .
\end{eqnarray*}
For $\alpha=0$ the problem is clearly integrable by using the
trivial set of first integrals $\{\tH,\psi_1,\psi_3\}$ in
involution (remarks~\ref{rem02} and \ref{rem03}). In the $\alpha
\ne 0$ case, one has the invariance-generators $\left\{
X_{2}=0,\Psi_{2}=0,T=2\,\lambda^{1}t+\lambda^{3},
\Psi_{1}=-\lambda^{1}\psi_{1},\right.
U_{1}=-\lambda^{1}u_{1},\Psi_{3}=-\lambda^{1}\psi_{3},
U_{2}=-\lambda^{1}u_{2}, X_{3}=\lambda^{1}x_{3}+\lambda^{2},
X_{1}=\lambda^{1}(\alpha^{-1}+ \left. x_{1}) \right\}$, which,
after solving the algebraic system of Theorem~\ref{main_teo} gives
the following set of effective first integrals
$$\left\{F_1=\tH,F_2=\left(\frac{1}{\alpha}+x_1\right)\psi_1
+x_3\psi_3-2t\tH,F_3=\psi_3\right\}.$$ This example shows that the
method not only generate mappings $F(x,\psi)$ verifying
$\{F,\tH\}=0$, but also, satisfying the more relaxed condition
$\{F,\tH\}=c\tH$ for some $c\in\bkR$ (see Remark~\ref{rem01}). We
observe that a nonautonomous first integral ($F_2$) is required to
prove integrability of the problem, in spite of the fact that the
problem is autonomous.
\end{example}


\section{The sub-Riemannian nilpotent case $(2,3,5)$}
\label{sec:SR235}

The sub-Riemannian (SR) problem concerns to characterize geodesics
in some $n$ dimensional SR-manifold $M$, i.e. to find absolutely
continuous curves $t\mapsto q(t)\in M$, $0\leq t\leq T$,
minimizing length
$$   l(q)=\int_0^T \,<\dot{q}(t),\dot{q}(t)>^{\frac{1}{2}}\,dt $$
such that $\dot{q}(t)\in \Delta(q(t))\backslash\{0\}$ a.e. $t$,
where $\Delta$ is a constant rank $m\leq n$ distribution with a
(degenerate) Riemannian metric $g$ on $\Delta$, and
$<\cdot,\cdot>$ is the scalar product induced by $g$. The
SR-problem can be locally formulated as an optimal control problem
\cite{bonnard}: let $(U,q)$ be a chart on which $\Delta$ is
generated by an orthogonal basis $\{X_1,\dots,X_m\}$, then the
SR-problem $(U,\Delta,g)$ is equivalent to
$$ \frac{1}{2} \int_0^T \left(\sum_{i=1}^m u_i^2(t)\right)\,dt\longrightarrow
\min\, , \quad \dot{q}(t)=\sum_{i=1}^m u_i(t)\,X_i(q(t)). $$ The
nilpotent case $(2,3,5)$ is the instance where $m=2$, $n=5$, and
the Lie algebra generated by $X_1$ and $X_2$ is a complete
nilpotent Lie algebra of nildegree $3$. On other words, if
$[\cdot,\cdot]$ denotes the Lie bracket of vector fields,
$X_3=[X_1,X_2]$, $X_4=[X_1,X_3]$ and $X_5=[X_2,X_5]$, the
SR-problem is nilpotent of type $(2,3,5)$ if
$X_i(0)=c_i\frac{\partial}{\partial x_i}$ for some $c_i\neq
0\in\bkR$ and $i\in\{1,\dots,5\}$; which gives
$\dim(\{X_1,X_2\})=2$, $\dim(\{X_1,X_2,X_3\})=3$ and
$\dim(\{X_1,X_2,X_3, X_4,X_5\})=5$. The true Hamiltonian is then
given by $\tH(q,\psi)=\frac{1}{2}\sum_{i=1}^m h_i^2(t)$, where
$h_i(t)=<\psi(t),X_i(t)>$ for $i\in\{1,\dots,m\}$, and the
Poincar\'e system is a system of equations on $T^\ast U$, given by
completing the set $\{X_1,\dots,X_m\}$ to form a smooth basis of
$TU$. Such vector fields are obtained by extending
$h_i(t)=<\psi(t),X_i(t)>$ to $i\in\{1,\dots,n\}$ and computing
$\dot{h}_i=\sum_{i=1}^m \{h_i,h_j\}\,h_j$, where $\{\cdot,\cdot\}$
is the n bracket.

Y. Sachkov proved in \cite{sachkov} that the optimal control
associated with the sub-Riemannian nilpotent case $(2,3,5)$ is
solvable, by obtaining three first integrals and using them to
reduce the Hamiltonian system to the differential equation
$\ddot{\theta}(t)=c_1\cos(\theta(t))+c_2\sin(\theta(t))$ for
$c_1,c_2,\theta(t)\in\bkR$. Hence, the solvability is obtained
from previous works showing that such differential equation is
integrable by quadratures using Jacobian Elliptic Functions. The
autonomous Hamiltonian $\tH$ is one of the first integrals used,
and the other two $h_4$ and $h_5$ are obtained as a direct
consequence of the fact that the Lie algebra is nilpotent, since
$\{h_i,h_j\}=<\psi,[X_i,X_j](q)>$ imply that $\dot{h}_4(t)=0$ and
$\dot{h}_5(t)=0$ along solutions.

With the method presented in this work, we can obtain enough first
integrals to directly prove the solvability of the problem. In
fact, we will consider a more general problem, parameterized by
constants $\alpha,\beta\in\{0,1\}$. Consider the SR nilpotent case
$(2,3,3+\alpha+\beta)$ with local generators for $\Delta$ given by
$$ X_1=\frac{\partial}{\partial x_1}
\:\:\mbox{ and }\:\:
 X_2=\frac{\partial}{\partial x_2}+x_1\frac{\partial}{\partial x_3}
+\frac{\alpha}{2}{x_1}^2\frac{\partial}{\partial x_4}+\beta x_1
x_2\frac{\partial}{\partial x_5}.$$ The distribution $\Delta$ is
known as the nilpotent Heisenberg distribution ($\alpha=0$ and
$\beta=0$), the nilpotent Engel distribution ($\alpha=1$ and
$\beta=0$), or the Cartan distribution ($\alpha=1$ and $\beta=1$).
The Pontryagin maximum principle gives the true Hamiltonian
$$\tH_{\alpha,\beta}=\frac{1}{2}\left[\psi_{1}^2+ (\psi_{2}
+x_{1}\psi_{3}+\frac{\alpha}{2}x_{1}^{2}\psi_{4}
+\beta x_{1}x_{2}\psi_{5} )^2\right].$$

This class of problems admit the following set of generators
(\textrm{cf.} \cite{GouveiaTorresRocha}):
$$\left\{ \Psi_{{5}}=-\frac{3}{4}\lambda^{{1}}\psi_{{5}},\Psi_{{1}}
=-\frac{1}{2}\lambda^{{1}}\psi_{{
1}},\Psi_{{2}}=-\frac{1}{2}\psi_{{2}}\lambda^{{1}},\Psi_{{3}}=-\lambda^{{1}}\psi_{{
3}}-\lambda^{{2}}\psi_{{5}}, \right.$$ $$\left.
T=\lambda^{{1}}t+\lambda^{{4}}, \Psi_{{4}}
=-\frac{3}{2}\lambda^{{1}}\psi_{{4}},X_{{1}}=\frac{1}{2}\lambda^{{1}}x_{{1}},
U_{{1}}=-\frac{1}{2}\lambda^{{1}}u_{{1}},X_{{2}}
=\frac{1}{2}\lambda^{{1}}x_{{2}}+\frac{1}{\beta}\lambda^{{2}},
\right.$$ $$\left. U_{{2}}=\frac{1}{2}\lambda^{{1}}u_{{2}},
X_{{5}}=\lambda^{{2}}x_{{3}}
+\frac{3}{2}\lambda^{{1}}x_{{5}}+\lambda^{{3}},X_{{4}}
=\frac{3}{2}\lambda^{{1}}x_{{4}}+\lambda^{{5}},X_{{3}}
=\lambda^{{1}}x_{{3}}+\lambda^{{6}} \right\}.$$ Computing
$\digamma(\lambda)$ on Lemma~\ref{lem01} and finding solutions for
Theorem~\ref{main_teo}, we have the effective first integrals (not
in involution)
$$\{\tH_{\alpha,\beta},\,\psi_{{2}}
+\beta \psi_{{5}} x_{{3}},\,\psi_{{3}},\,\alpha\psi_{{4}},\,\beta\psi_{{5}}\}.$$

\begin{lemma}
The sub-Riemannian nilpotent cases $(2,3)$
($\alpha=0$ and $\beta=0$), $(2,3,4)$ ($\alpha=1$ and $\beta=0$),
and $(2,3,5)$ ($\alpha=1$ and $\beta=1$), are integrable by quadratures.
\end{lemma}

It is not difficult to find a first integral $F$ whereas the set
$\{\tH,\,F,\,\psi_{{3}},\,\psi_{{4}},\,\psi_{{5}}\}$ is
involutive. Such first integral should satisfy the relations
$\{\tH,F\}=0$ and $\{F,\psi_i\}=0$ for $i\in\{3,4,5\}$. If we
consider \emph{a priori} that $F$ does not depend on $x_3$, $x_4$
or $x_5$, then last condition is trivially verified. Hence, it
just remains to solve $\{\tH,F(x_1,x_2,\psi_1,\dots,\psi_5)\}=0$,
which by a direct calculation gives the first integral
$$ F=-\psi_1\psi_5+\psi_2 \psi_4-(\psi_3+\frac{1}{2}\psi_5
x_2)x_2\psi_5.$$ Therefore, the solutions of the sub-Riemannian
nilpotent case $(2,3,5)$ are Liouvillian, using the first
integrals $\{\tH,\,F,\,\psi_{{3}},\,\psi_{{4}},\,\psi_{{5}}\}$ in
involution.

\bigskip

Although the present method can be applied to other hard problems,
such as the sub-Riemannian nilpotent cases $(2,3,5,8)$ or
$(2,3,5,8,14)$, for which the solvability is still unknown,
because of its complexity (8 and 14 effective first integrals are
needed, respectively), their study is left for a forthcoming
publication. Here we just notice that the Maple package described
in \cite{GouveiaTorresRocha} is unable to find a sufficient rich
family of first integrals for the case $(2,3,5,8,14)$. Therefore,
they need to be found by other theoretical procedure.


\section*{Acknowledgments}

This work was partially supported by the Portuguese Foundation for
Science and Technology (FCT), through the Control Theory Group
(cotg) of the Centre for Research in Optimization and Control
(CEOC).



\end{document}